\documentclass[twoside]{irmaems}
\usepackage{amssymb} 
\usepackage{amsmath} 
\usepackage{latexsym}

\usepackage{caption, subcaption, amsfonts}
\usepackage{ color, amsmath, amscd, graphicx,latexsym, hyperref}

\newcommand{\sherh}[1]{\fboxsep=0pt\setlength{\fboxrule}{1pt}
\begin{center}
   \fbox{\colorbox{green}{
         \begin{minipage}[t]{13cm}
            #1
         \end{minipage}
      }
   }
\end{center}}
\renewcommand{\sherh}[1]{}

\renewcommand\theequation{\arabic{section}.\arabic{equation}}

\theoremstyle{definition} 

 \newtheorem{definition}{Definition}[section]
 
 \newtheorem{example}[definition]{Example}

\theoremstyle{plain}      

 \newtheorem{theorem}[definition]{Theorem}

\newcommand{\mcg}{\operatorname{MCGroupoid}}
\renewcommand{\mod}{\operatorname{Mod}}

\newcommand{\ITEM}{$\blacktriangleright$}

\newcommand{\TM}{Teichm\"{u}ller}

\newcommand{\hh}{{\mathbb H}}
\newcommand{\CC}{\mathbb{C}}
\newcommand{\ZZ}{\mathbb{Z}}
\newcommand{\QQ}{\mathbb{Q}}
\newcommand{\RR}{\mathbb{R}}
\newcommand{\HH}{\mathbb{H}}
\newcommand{\PP}{\mathbb{P}}
\newcommand{\mm}{\mathcal{M}}
\newcommand{\DD}{\mathbb{D}}
\newcommand{\F}{\mathcal F}
\newcommand{\lra}{\longrightarrow}
\newcommand{\mat}[2]{\left (\begin{array}{c c}
					#1 \\ #2
	          \end{array}
\right)}
\newcommand{\sm}{\setminus}
\newcommand{\cov}{{\rm {\bf Cov}}}
\newcommand{\fcov}{{\rm {\bf FCov}}}
\newcommand{\fsub}{{\rm {\bf FSub}}}
\newcommand{\fgicov}{{\rm {\bf FGICov}}}

\newcommand{\sub}{{\rm {\bf Sub}}}
\newcommand{\psl}{\mathrm{PSL}(2,\ZZ)}
\newcommand{\gal}{G_{\mathbb Q}}
\def\modorb{\mbox{$\circ\hspace{-1.5mm}-\!\!\!-\hspace{-1.5mm}\bullet$}}
\newcommand{\wt}[1]{\widetilde{#1}}

\begin{document}

\title{A panaroma of the fundamental group of the modular orbifold}

\author{A. Muhammed Uluda\u{g}\thanks{
Work supported by a Galatasaray University Research Grant and the grant TUBITAK-114R073} and  Ayberk Zeytin\thanks{Work
supported by T\"{U}B\.{I}TAK Career Grant No.~113R017 and 114R073}}

\address{
Department of Mathematics, Galatasaray University\\
\c{C}{\i}ra\u{g}an Cad. No. 26 Be\c{s}ikta\c{s} 34357 \.{I}stanbul, Turkey\\
email:\,\tt{muludag@gsu.edu.tr}
\\[4pt]
Department of Mathematics, Galatasaray University\\
\c{C}{\i}ra\u{g}an Cad. No. 26 Be\c{s}ikta\c{s} 34357 \.{I}stanbul, Turkey\\
email:\,\tt{azeytin@gsu.edu.tr}
}

\maketitle

\begin{abstract} 	We give an overview of the category of subgroups of the modular group, incorporating both the \emph{tame} part, i.e. finite index subgroups, and the \emph{non-tame} part, i.e. the rest. We also discuss arithmetic related questions which exist in both the tame part (via Belyi's theorem) and the non-tame part.
\end{abstract}



\tableofcontents   

\section{Introduction}\label{sec:intro}
The fundamental group of the modular orbifold, $\mm$, is the modular group, $\psl$. Due to the general topological facts, there is a one to one correspondence between subgroups of the modular group and covers of $\mm$. Denote by ${\rm {\bf Sub}}(\psl)$ the category of all subgroups of the modular group with inclusions as morphisms and by ${\rm {\bf FSub} }(\psl)$ the category of all finite-index subgroups of $\psl$, with inclusions as morphisms. These two categories are equivalent to the base-pointed covering category $\cov^{*}(\mm)$ (resp. $\fcov^*(\mm)$) of (resp. finite) covers of the modular orbifold. Denote by $\cov(\mm)$ and $\fcov(\mm)$ the corresponding covering categories obtained by forgetting the base points. Our aim here is to give the reader a panorama of these categories.

The categories $\cov(\mm)$ and $\cov^{*}(\mm)$ and $\sub(\psl)$ contains uncountably many objects. On the other hand, objects in the categories $\fcov(\mm)$ and $\fsub(\mm)$ can be represented by modular graphs, i.e. a certain class of finite bipartite ribbon graphs or dessins d'enfants, usually thought of as graphs embedded on oriented topological surfaces. By a celebrated theorem of Belyi, \cite{belyi}, these graphs admit a Galois action, i.e. an action of the absolute Galois group, $\gal$. This action is faithful. 
The idea of studying $\gal$ through this interaction between $\gal$ and the {modular group} $\psl$ 
is the main motivation of the ``dessins" program, \cite{grothendieck-esquisse,grothendieck/marche}. However, Deligne put already in 1989: ``Grothendieck and his students developed a combinatorial description (``maps") of finite covers of $\PP^1_\CC$... This did not help to understand the Galois action. One has just a few examples of non-solvable covers whose Galois conjugates have been computed", \cite{le/groupe/fondamental/deligne}. 


Confirming Deligne's convictions,  it seems that the picture that emerged after almost 40 years of research is far from the original expectations expressed very lively in \cite{grothendieck-esquisse} and \cite{grothendieck/marche}. It turns out that understanding this Galois action on all of ${\rm {\bf FCov} }( \mathcal M)$ appears to be a hopeless task. One turns attention to some manageable substructures of it.
In the previous chapter, a certain sub-structure (the family of hypergeometric covers) of $\fcov(\mm)$ is suggested to be studied from the arithmetic point of view. In the current chapter, we suggest to go beyond $\fcov(\mm)$, and study some infinite covers in $\cov(\mm)$. It turns out that although of completely different nature, arithmetic is still present for infinite index subgroups, see for instance \cite{UZD} and \cite{reduction}. 
\section{Modular graphs}\label{sec:modular/graphs}
Let $\mathbb H$ be the upper half plane with the standard action of the modular group  $\psl$; i.e. the group of invertible $2\times2$ matrices with integer coefficients. This action sends $z\in \mathbb H$ to 
$$\mat{p&q}{r&s} \cdot z \mapsto \frac{pz + q}{rz + s}.$$
Its quotient space  is called the {\it modular orbifold}\index{modular orbifold} and denoted $\mathcal M$:
$$
\mathcal M=\psl\backslash\mathbb H
$$
The modular orbifold has two orbifold points, $i$ with stabilizer {\small $\langle S =  \mat{0 & -1}{1 & 0}\rangle \cong \ZZ \!/\! 2\ZZ$} and 
$\exp{2 \pi i/3}$ with stabilizer {\small $\langle L = \mat{1 & -1}{1 & 0}\rangle \cong \ZZ \!/\! 3\ZZ$.} Moreover, its orbifold fundamental group is
$$
\pi_1(\mathcal M)=\psl\simeq \ZZ\!/\!2\ZZ \ast \ZZ\!/\!3\ZZ.
$$ 
By the usual correspondence from topology, the subgroups of the fundamental group classify the coverings of $\mathcal M$. 
The underlying space of $\mathcal M$ is the projective line
$\PP^1$, the two orbifold points with $\ZZ\!/\!2\ZZ$ and $\ZZ\!/\!3\ZZ$-inertia may be assumed to be respectively at $0$ and $1$. The translation $\mat{1&1}{0&1}$ is responsible for the puncture of the quotient which may be taken to be the point at $\infty$.  There is a degree-six Galois covering of the modular orbifold by the projective line punctured at three points; $\PP^1\sm\{0,1,\infty\}$. Therefore every covering of $\PP^1\backslash\{0,1,\infty\}$ is also a covering of $\mm$. Note that every covering of $\mathcal M$ carries a canonical hyperbolic metric induced from the Poincar\'{e} metric on the upper half plane.

\begin{figure}[h!]
	\centering
	\includegraphics{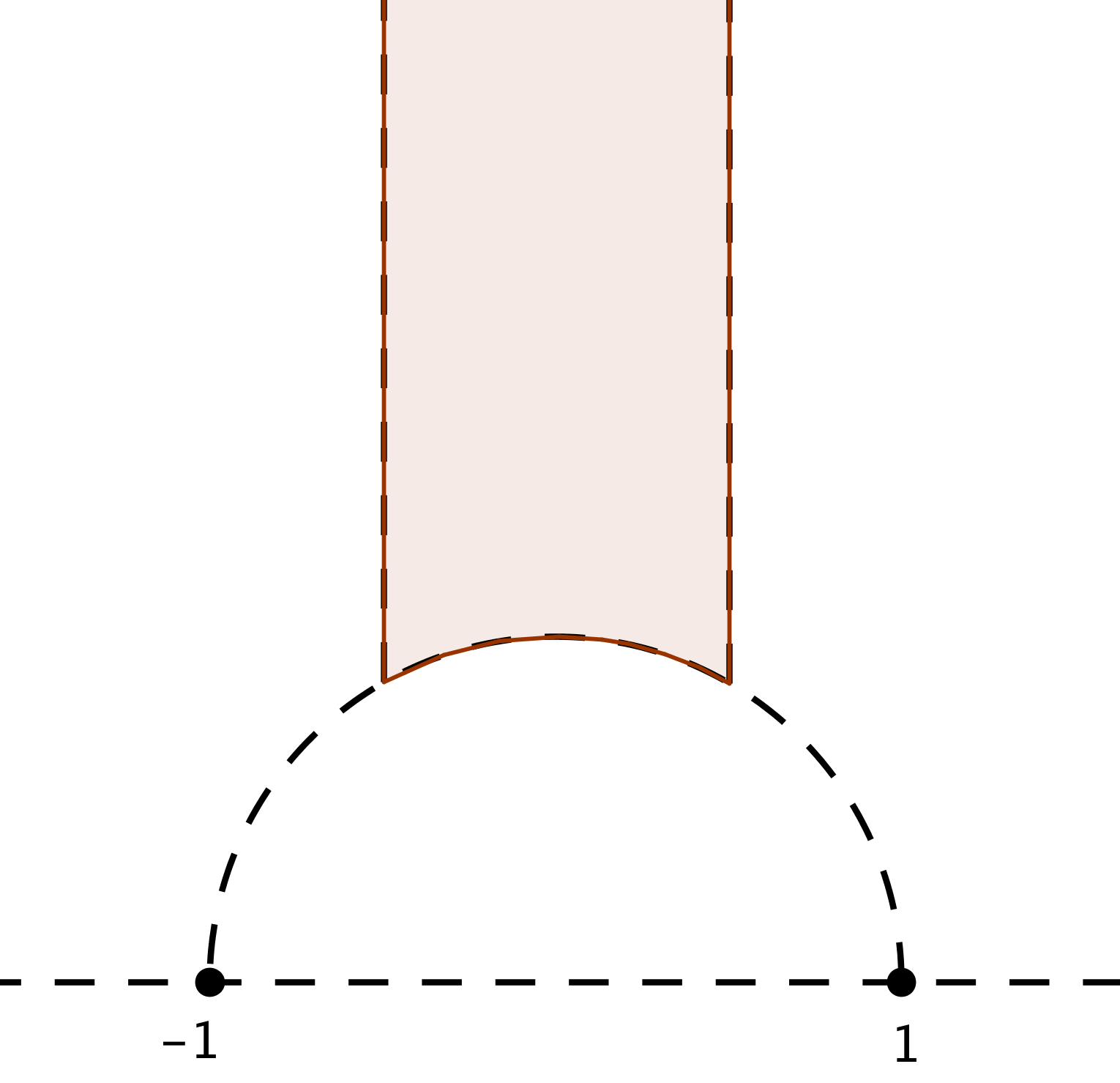}
	\caption{The standard fundamental domain for the action of $\psl$ on the upper half plane.}
	\label{fig:fundamental/region}
\end{figure}

Finite covers of the modular orbifold admit a combinatorial description by certain graphs (Linienzuges of Klein, \cite{klein/linienzuges}), as follows. Consider the arc connecting the two elliptic points on the boundary of the standard fundamental domain of the $\psl$ action on $\mathbb H$, see Figure~\ref{fig:fundamental/region}. 

Then the $\psl$-orbit of this arc is a tree $\mathcal F$, which we call the {\it Farey tree}\index{Farey tree}. Now by its very definition, this tree admits a $\psl$-action, and a quotient graph by a subgroup of finite or infinite index $G<\psl$ has been named by us a {\it modular graph}\footnote{Modular graphs are very similar to trivalent ribbon graphs. This new baptisation is our modest contribution to the very rich and diversified terminology concerning these combinatorial objects.}\index{modular graph} \cite{UZD}. Thus the modular graph $G\backslash\mathcal F$ sits inside the curve $G\backslash\HH$ in a standard way. In particular, the quotient \emph{orbi-graph} $\psl \backslash\mathcal F$ is an arc connecting the two orbifold points of the modular orbifold $\psl \backslash\mathbb H$. We call this the {\it modular arc}\index{modular arc} and denote it by $\modorb$. 

From another perspective, modular graphs are coverings (in the orbifold sense) of the modular arc $\modorb$. The degree-$n$ covers of $\modorb$ are obtained by gluing $n$ copies of $\modorb$ at the endpoints, to obtain a connected graph, so that 

\begin{enumerate}
	\item endpoints of different types never meet, 
	\item type-$\circ$ vertices are of degree either 1 or 2, and type-$\bullet$ vertices are of degree either 1 or 3, and 
	\item for each vertex of degree 3, a cyclic order of the edges meeting at that vertex is given.
\end{enumerate}

A systematic study of these graphs in connection with modular curves was made by Kulkarni, under the name cyclic trivalent graphs and tree diagrams. The Farey tree $\mathcal F$ is the universal cover of $\modorb$. 

Perhaps the simplest modular graphs that one may draw this way (as coverings of the modular interval) are finite planar trivalent trees, (with two types of terminal vertices and a puncture at infinity). Hence we may consider finite planar trivalent trees as natural objects in $\fsub(\psl)$. Their Belyi maps are polynomials generalizing the Chebyshev polynomials.

One may imagine widening the puncture of the modular orbifold $\mathcal M$, until the modular orbifold ``deformation retracts'' to $\modorb$ inducing an isomorphism of fundamental groups. Let us observe what happens to the coverings during this deformation retraction: at time $0$ we start with punctured surfaces, at time $1-\epsilon$ we see ribbon graphs and at time $1$ we get the cyclic trivalent graphs. Meanwhile the universal cover of $\mathcal M$, i.e. the upper half plane, deformation retracts to the Farey tree. One may imagine that the category ${\rm {\bf FCov} } (\mathcal M)$ itself is deformation retracting to ${\rm {\bf FCov} } (\mathcal \modorb$), see Figure~\ref{fig:defret}. Naming the inverse limit of  ${\rm {\bf FCov} } ^*(\modorb)$ as the {\it ribbon solenoid}, this operation gives a deformation retract of the punctured solenoid to the ribbon solenoid, which is a compact space. So we obtain:
\begin{figure}[h!]
	\begin{center}
		\includegraphics[scale=1]{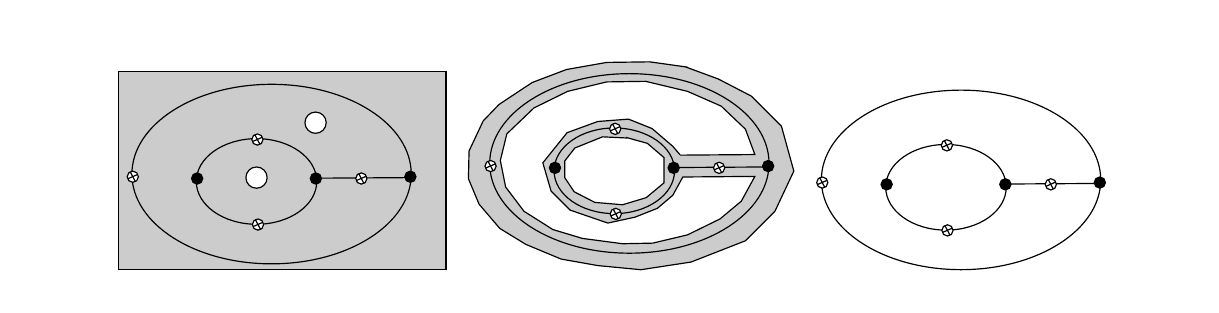}
	\end{center}
	\caption{Deformation retract of the ambient surface to a modular graph.}
	\label{fig:defret}
\end{figure}
\begin{theorem}
	The pointed categories ${\rm {\bf FCov}^* }(\mathcal M)$ and ${\rm {\bf FCov} } ^*(\modorb)$ (= pointed finite modular graphs) are equivalent via deformation retraction, and they are both arrow-reversing equivalent to  ${\rm {\bf FSub} }(\psl)$ via the fundamental group functor.
\end{theorem}
The category of modular graphs admits a {\it standard} realization by piecewise geodesic and analytic arcs of ${\rm {\bf FCov}^* }(\mathcal M)$, via
$$
G\backslash \mathcal F \subset G\backslash \mathbb H.
$$
Since there are countably many finite modular graphs, these categories admit countably many objects. The formula $\#\{\mbox{vertices}\} - \#\{\mbox{edges}\} + \#\{\mbox{faces}\}$\footnote{A face on $G \backslash\mathcal{F}$ is a finite left-turn closed path.} determines the Euler characteristic and hence the genus of the curve, whenever the graph is finite. 

In fact, we haven't made use of any finiteness property in the proof, and therefore conclude more generally that the categories ${\rm {\bf Cov}^* }(\mathcal M)$,  ${\rm {\bf Cov} } ^*(\modorb)$ and ${\rm {\bf Sub} }\, \psl$ are equivalent as well. There are uncountably many infinite modular graphs, so there are as many infinite coverings of the modular orbifold, equivalently, as many subgroups of the modular graphs. 

If one \emph{erases} vertices of degree two from a modular graph then the remaining graph becomes trivalent, i.e. every vertex is of degree $3$. 
Hence, trivalent ribbon graphs \cite{penner/decorated/TM/theory}  are inside ${\rm {\bf FCov} }(\modorb)$ in a natural way. There are more modular graphs than  trivalent ribbon graphs since the former are allowed to have terminal edges abutting at orbifold points of both types. On the other hand, ``{a trivalent ribbon graph (or an ideal triangulation) with a distinguished oriented edge}" encountered frequently in {\TM} theory (both in the finite part, see \cite{penner/decorated/TM/theory, ptolemy} and the infinite part, see \cite{universal/constructions/penner,imbert/isomorphisme/du/T/avec/P}) is same as the data ``{a modular graph with the choice of a base edge}". Note that the mapping class group acts freely on the set of embedded  trivalent ribbon graphs, and the groupoid associated to this group action admits a nice presentation with modular graphs as objects and graph flips as generator of morphisms, see \cite{penner/decorated/TM/theory}.

Every triangulation of a surface has a dual graph which is an embedded trivalent ribbon graph. If we include degenerate triangulations this correspondence extends to modular graphs. This correspondence is bijective between the orientation-preserving homeomorphism classes of triangulations and trivalent ribbon graphs.

Conversely, any modular graph determines a covering of $\PP^1$ branched only at $0,1$ and $\infty$. Indeed, given a modular graph one chooses a base edge $e_{o}$ on the graph. Then, every finite loop\footnote{Every terminal vertex is considered as an orbifold point. Hence the path from $e_{o}$ to a terminal vertex and back is a generator of finite order of the subgroup.} based at $e_{o}$ on the graph determines an element, $W\in \psl$ in the following way: starting with the empty word $W$ while tracing the loop if a degree two vertex is visited one adds an $S$ to $W$ and if a degree three vertex is visited then one writes $L$ (resp. $L^{2}$) if the next edge on the loop is on the left (resp. right) of the initial edge. Then, all such loops generate a subgroup, say $G<\psl$, and as a result of the construction explained above, the orbifold $\psl \backslash\mathbb{H}$ deformation retracts onto the modular graph. In fact, the Riemann surface $G \backslash\HH$ can be defined over $\overline{\QQ}$. 

\subsection{Modular tile and the holy triality.}
\begin{figure}[h!]
	\centering
	\begin{subfigure}{4cm}
		\includegraphics[scale=0.18]{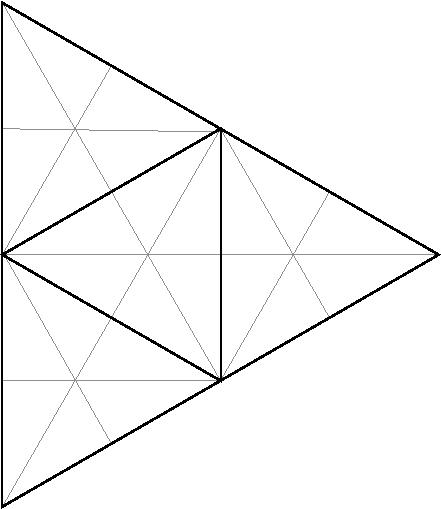}		
		\caption{A triangulation}
	\end{subfigure} \quad
	\begin{subfigure}{4cm}
		\includegraphics[scale=0.18]{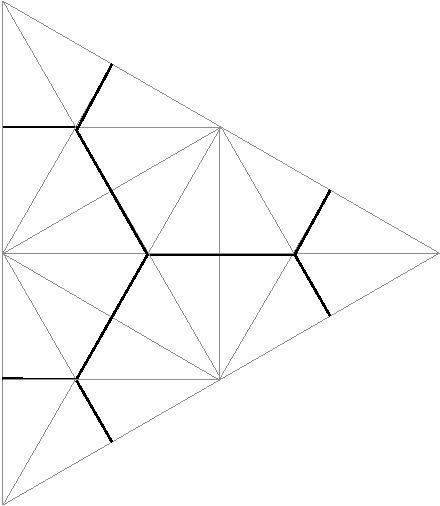}		
		\caption{its modular graph}
	\end{subfigure} \quad
	\begin{subfigure}{4cm}
		\includegraphics[scale=0.18]{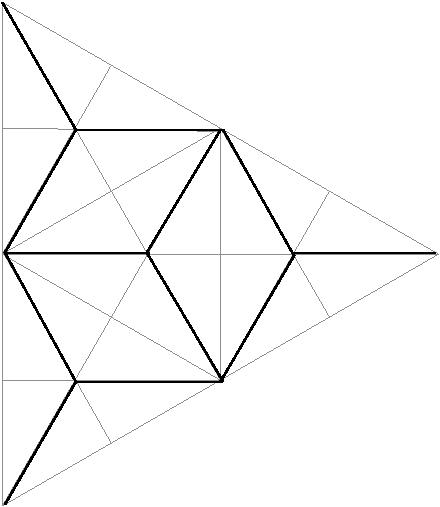}
		\caption{and its lozenge}
	\end{subfigure}
	\caption{Triality of graphs. Small tiles are the copies of the modular tile, rendered here under the flat metric. Vertices of the triangulation are punctures.}
	\label{fig:triality}
\end{figure}

Every modular graph $G\backslash \mathcal F$ has a standard piecewise analytic realization on the Riemann surface $G\backslash {\mathbb H}$ with edges being geodesic segments with respect to the hyperbolic metric induced from the upper half plane. Equivalently, these edges are lifts of the modular arc by the canonical projection $G  \backslash {\mathbb H} \longrightarrow \psl\backslash {\mathbb H}$. If instead we lift the geodesic arc connecting the $\ZZ\!/\!3\ZZ$-elliptic point to the cusp to the surface $\psl  \backslash {\mathbb H}$, then we obtain another graph, which is the associated {\it ideal triangulation}\index{ideal triangulation}, i.e. a triangulation of the surface whose vertices are exactly at the punctures. Lifting the remaining geodesic arc gives rise to yet another type of graph, called a {\it lozenge tiling}\index{lozenge tiling}. So there is a triality, not just duality, of these graphs, see Figure~\ref{fig:triality}.

In fact, for some purposes it is more natural to study the ${\rm \bf FSub}$ of the extended modular group $\mathrm{PGL}_2(\ZZ)$. This
is the fundamental group of the quotient  $\mathcal M^*$ of the modular orbifold $\mathcal M$ by the complex conjugation. 
We shall call $\mathcal M^*$ the {\it modular tile}. 
Elements of ${\rm \bf FCov}(\mathcal M^*)$ are surfaces tiled by copies of $\mathcal M^*$ and they 
include bordered surfaces with some corners and punctures, possibly lying on the boundary.
Every element of ${\rm \bf FCov}(\mathcal M^*)$ has a double which is in ${\rm \bf FCov}(\mathcal M)$.
Any polygon, which is triangulated by its diagonals is an element of ${\rm \bf FCov}(\mathcal M^*)$.
Lifts of the each one of the three edges of the modular tile to its coverings 
gives rise to modular graphs, triangulations, and lozenges.

\subsection{Special modular graphs}
Now that we have introduced modular graphs, the next task is to classify them. There are several schemes that come to mind, the first being with respect to genus and the number of punctures/orbifold points of the ambient surface $G\backslash \mathbb H$. It seems natural to denote this slice of the covering category by 
$$
{\rm \bf FCov}_{d, g, n, n_\circ, n_\bullet}(\mathcal M).
$$
Here $d$ is the number of edges of the modular graph (= index of a corresponding subgroup), $g$ is the genus of the ambient Riemann surface, $n$ is the number of punctures, $n_{\circ}$ is the number of orbifold points of type $'circ$, $n_{\bullet}$ is the number of orbifold points of type $\bullet$. Note that these invariants are not independent and they can be defined in a purely combinatorial manner. The genus-0 case is the case of modular graphs on the sphere. The lowest value that the number of punctures $n$ can attain is 1. When the genus is 0, the case $n=1$ yields loop-free modular graphs, i.e. planar bipartite trees with vertices of type $\circ$ being of degree 1 or 2 and those of type $\bullet$ being of degree 1 or 3. We may call them ``{modular trees}''. If we erase vertices of degree 2, these are precisely the planar trivalent trees. 

Reasoning geometrically or combinatorially, one may imagine other ways to classify the modular graphs and operations to produce new graphs from old ones, \cite{wood/belyi/extending/maps}. One possibility is to take a triangulation, simultaneously subdivide each edge into $n$ segments and then simultaneously subdivide every triangle into $n^2$ smaller triangles by connecting the endpoints of these segments in the obvious manner. In the previous chapter it is shown that this operation is of some relevance in the theory.

\subsubsection{Congruence modular graphs.}
There is a well-known and important class of subgroups in the panorama of the modular group. 
These are defined as, for positive integer $N$
\begin{eqnarray*}
\Gamma_0(N)=\left\{ \mat{a & b}{c & d}\in \mathrm{SL}_2(\mathbb Z) \, |\, c\equiv 0 \mod N \right\}\qquad\qquad\qquad\quad\\
\Gamma_1(N)=\left\{ \mat{a & b}{c & d}\in \mathrm{SL}_2(\mathbb Z) \, |\, c\equiv 0, \quad a\equiv d\equiv 1 \mod N \right\}\quad\\
\Gamma(N)=\left\{ \mat{a & b}{c & d}\in \mathrm{SL}_2(\mathbb Z) \, |\, b\equiv c\equiv 0, \quad a\equiv d\equiv 1\mod N \right\}.
\end{eqnarray*}
$\Gamma(N)$ is called the {\it principal modular congruence group} and $\Gamma_0(N)$ and $\Gamma_1(N)$ are said to be {modular groups of Hecke type}. The same definitions apply to the projectivizations. A subgroup of the modular group containing a principal congruence subgroup is called a {\it congruence modular group}, and the corresponding coverings of the modular orbifold are called {\it modular curves}. For details 
(which occupy a considerable volume in the mathematical literature) about these subgroups we refer the reader to 
\cite{miyake/modular/forms}.

Principal congruence subgroups are the kernels of the reduction maps onto $\mathrm{SL}_2(\mathbb Z/N\mathbb Z)$ and therefore they are normal subgroups. Hence they correspond to Galois coverings of the modular orbifold and provide us with a family of modular graphs with a transitively acting monodromy group.

\subsection{Non-tame subgroups}
The finite part $\fcov(\mm)$ (or $\fsub(\psl)$) of $\cov(\mm)$ (or \linebreak $\sub(\psl)$) will be referred to as the \emph{tame} part of the category because of the fact that these give rise to \'{e}tale covers of the modular orbifold and conversely.  All the remaining objects in the category will be called \emph{non-tame}. 
As we noted above, the full category $\mathbf{Cov}(\mathcal M)$ admit uncountably many objects and what is beyond $\mathbf{FCov}(\mathcal M)$ appears to be a wild territory. However, there is a subcategory of $\mathbf{Cov}(\mathcal M)$ which consists of infinite covers of finite topology (i.e. with finitely many punctures, boundary components and handles). 
We denote this sub-category by $\mathbf{FGICov}(\mathcal M)$ (Finitely Generated Infinite Covers). Once again there is the pointed category $\mathbf{FGICov}^*(\mathcal M)$. As the name suggests, 
the corresponding subgroup category $\mathbf{FGISub}(\psl)$ consists of finitely generated subgroups of infinite index in the modular group. Unlike $\mathbf{FCov}(\mathcal M)$, it has no initial objects as the modular orbifold itself is not included in it. As a second remark, note that an infinite cover of an object of $\mathbf{FGICov}(\mathcal M)$ may or may not be inside $\mathbf{FGICov}(\mathcal M)$. On the other hand, this category is closed under finite covers. However, it is not true that any two objects have a common finite cover (i.e. the intersection of two finitely generated infinite subgroups is not always a finitely generated subgroup\footnote{For a concrete example,  take the two $\ZZ$-subgroups $\langle \gamma_1\rangle$ and $\langle \gamma_2\rangle$, where $\gamma_1$ and $ \gamma_2$ are two distinct primitive elements of infinite order.}), so one can not pass to the limit of $\mathbf{FGICov}(\mathcal M)^*$ to define some sort of profinite completion. 

We see that one may consider two versions of $\mathbf{FGICov}(\mathcal M)^*$, one with finite covers and the other with finite or infinite covers as morphisms. The idea is that if we restrict to finite covers then the remaining category will become entirely ``{arithmetic}". Taking an object in this restricted category and considering the subcategory of Galois coverings of this object, one may define the fundamental group of this object to be the limit of the inverse system of Galois groups of these normal covers.

\begin{figure}[h!]
	\centering
	\includegraphics[scale=1.7]{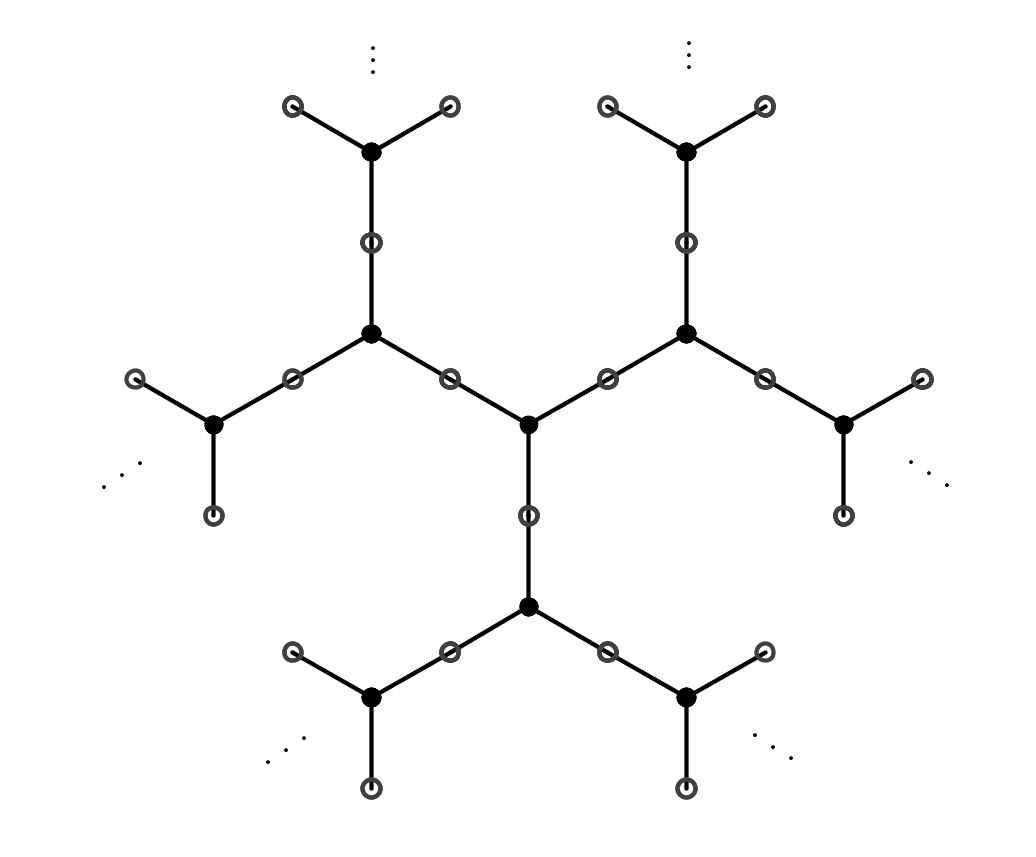}
	\caption{The bipartite Farey tree.}
	\label{fig:farey/tree}
\end{figure}

\begin{figure}[h!]
	\centering
	\includegraphics{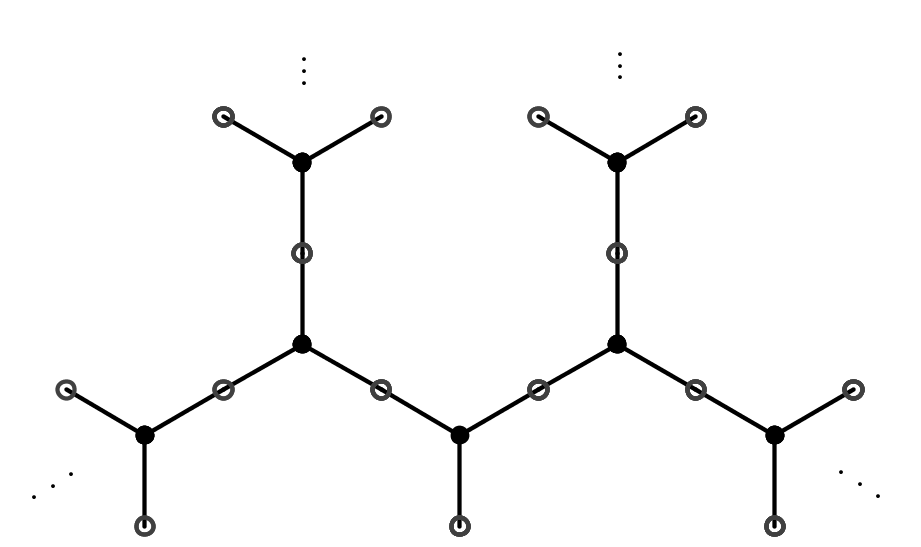}
	\quad
	\includegraphics{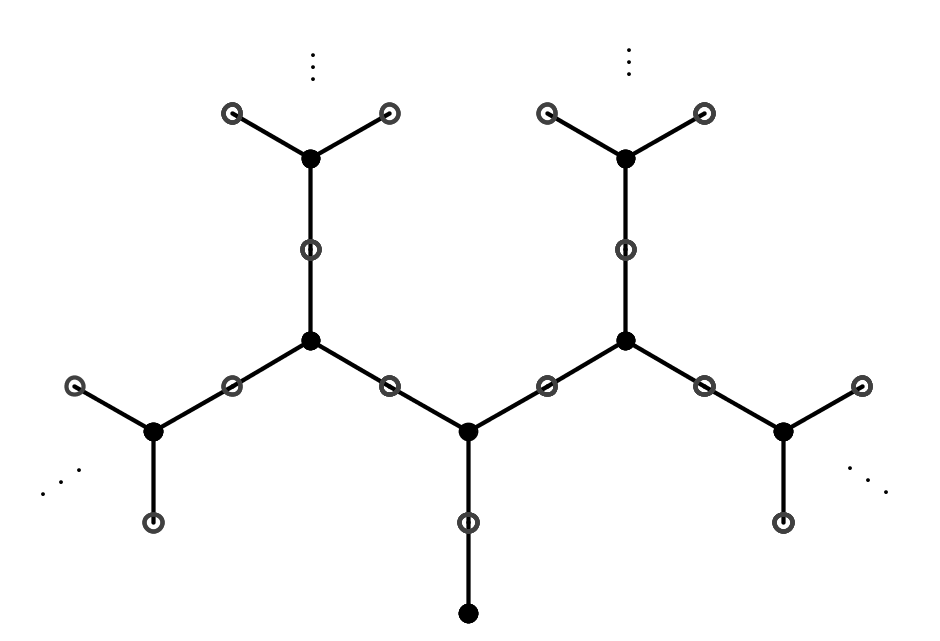}
	\caption{The modular graphs for $G = \{I,S\}$ and for $G = \{I,L,L^{2}\}$, respectively.}
	\label{fig:elliptic/carks}
\end{figure}

\begin{figure}[h!]
	\centering
	\includegraphics[scale=0.7]{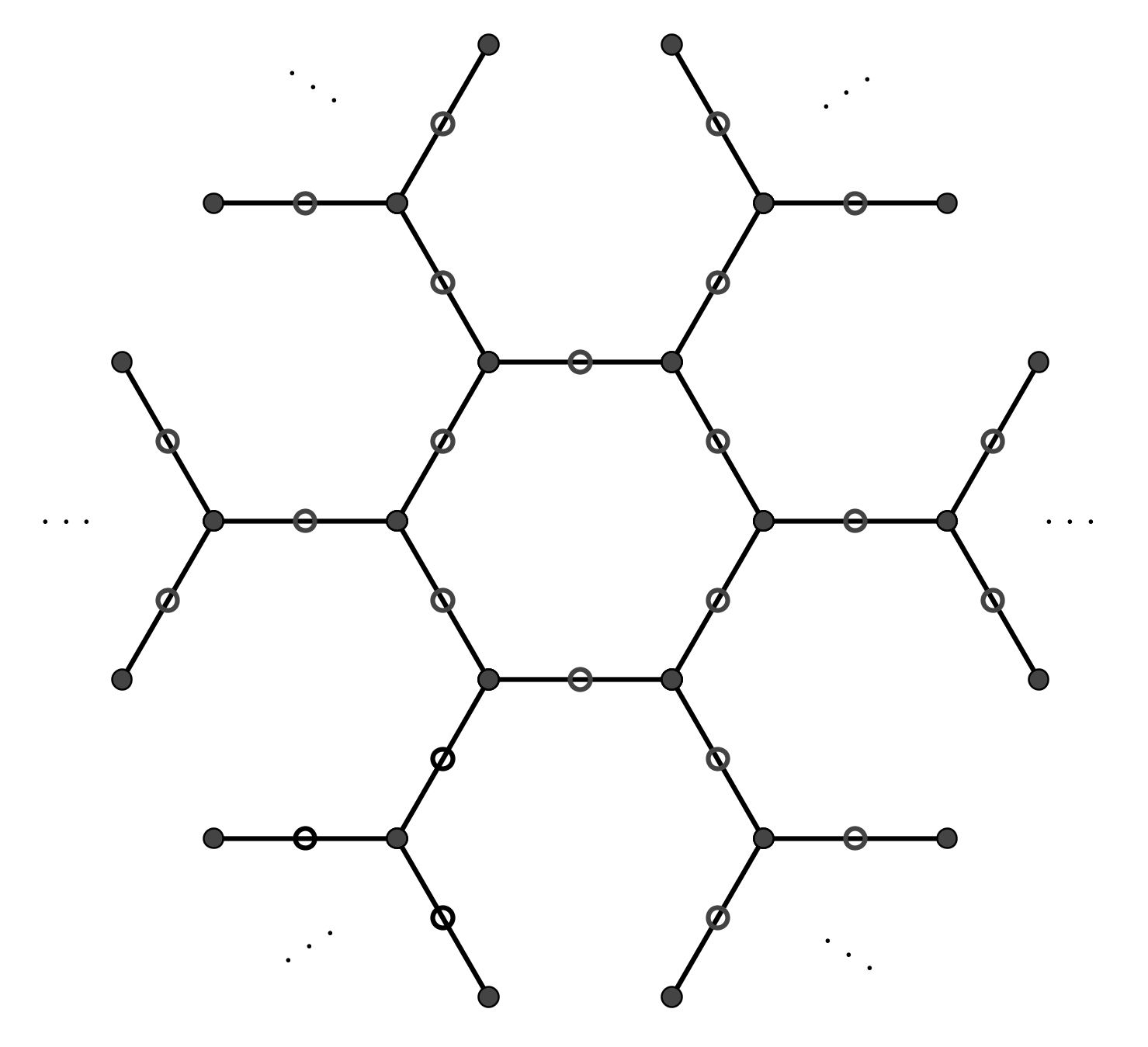}
	\caption{The modular graph corresponding to $G = \langle (LS)^{6} \rangle$.}
	\label{fig:parabolic/carks}
\end{figure}

\begin{figure}[h!]
	\centering
	\includegraphics[scale=0.7]{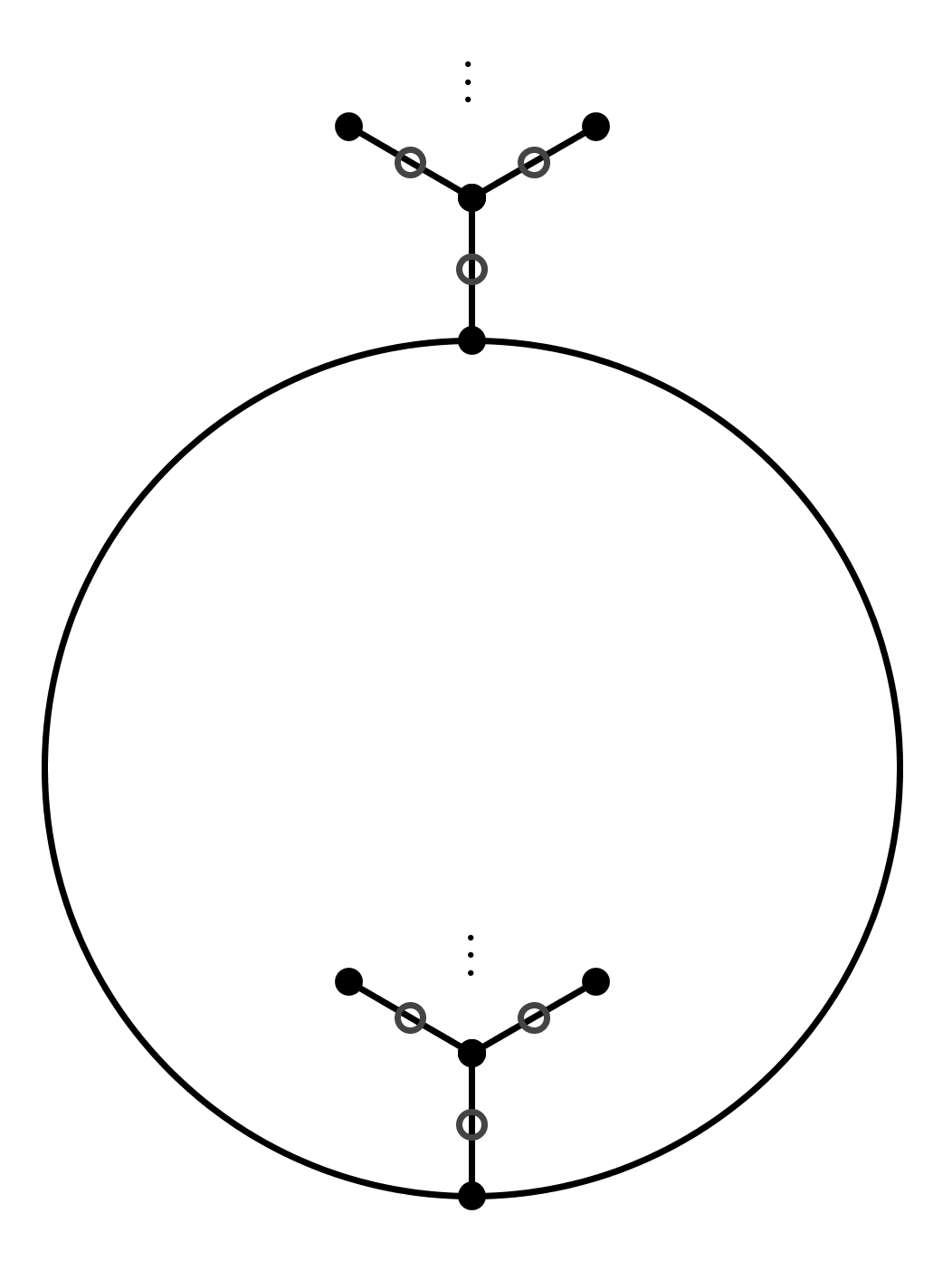}
	\caption{The modular graph corresponding to $G = \langle LSL^{2}S \rangle$.}
	\label{fig:hyperbolic/carks}
\end{figure}

\begin{figure}[h!]
	\centering
	\includegraphics[scale=0.5]{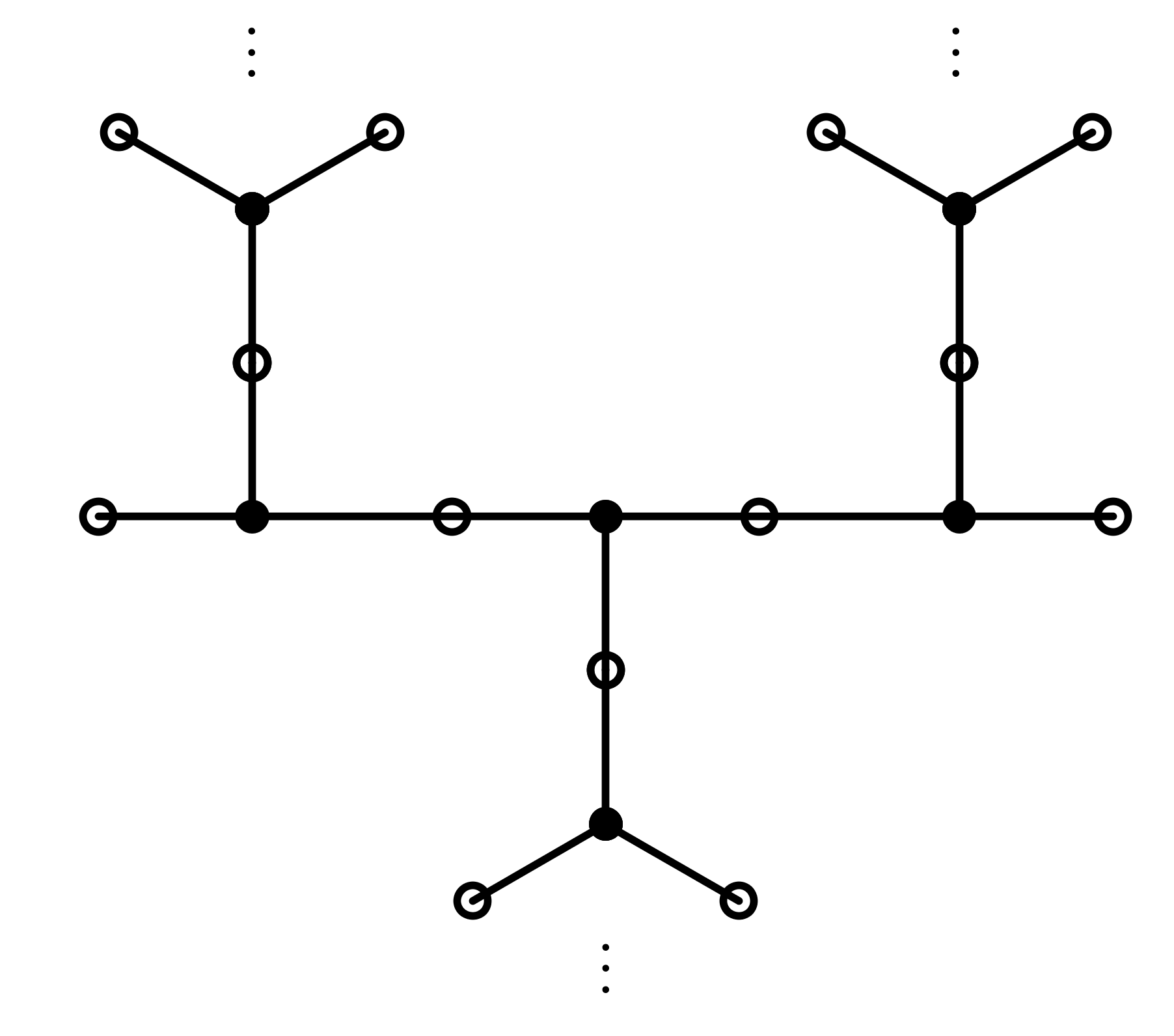}
	\caption{The modular graph of a dihedral subgroup}
	\label{fig:reciprocal}
\end{figure}

In order to demystify what have been said in the above paragraphs, let us consider some examples:

\begin{example}
	The very first example is the upper half plane itself, corresponding to the trivial subgroup, which is indeed finitely generated of infinite index. Observe that this is the final object of $\mathbf{FGICov}(\mathcal M)^*$ (if we admit infinite coverings inside the category). The corresponding graph object is the bipartite Farey tree itself, see Figure~\ref{fig:farey/tree}. Observe that the action of the modular group can also be described as an action on the bipartite Farey tree. This is very useful in considering graphs of arbitrary subgroups.
	\label{ex:trivial/subgroup}
\end{example}
\begin{example}[Finite Subgroups]
	As the next example let us consider the elliptic subgroup $\{1,S\}$. In this case, the graph has infinitely many vertices, and $S$ acts on the Farey tree by rotation about one of the degree $2$ vertices of the graph. Therefore $\Gamma _{G}$ is obtained by folding the Farey tree from the \emph{symmetry} vertex, see Figure~\ref{fig:elliptic/carks}. In a similar fashion the matrix $L$ acts on the Farey tree as a rotation around a degree $3$ vertex of degree $2\pi/3$, see Figure~\ref{fig:elliptic/carks} for the graph corresponding to $\{I,L,L^{2}\}$. This depicts the general case for subgroups generated by one element of finite order. Indeed, every such element is conjugate to either $S$ or $L$.
	\label{ex:finite/subgroups}
\end{example}
\begin{example}[Parabolic Subgroups]
	An element $W$ of $\psl$ is called parabolic if its absolute trace is $2$. As was the case for finite subgroups, any such element is conjugate either to $(LS)^{k}$ or $(L^{2}S)^{k}$ for some integer $k$. The corresponding graph has a unique loop (called the \emph{spine}) containing $2k$-many edges. There are $k$-many Farey tree components (called a Farey branch) attached to the degree $3$ vertices of the spine (there are $k$ such vertices!). Observe that the Riemann surface $\langle W \rangle \backslash \hh$ is a punctured disk. As was observed previously punctures on the Riemann surface are in one to one correspondence with finite left-turn closed paths on the graph. Therefore all Farey branches appear either on the left or on the right of the spine, depending on the orientation and they \emph{expand} in the direction of the \emph{outer boundary}, see Figure~\ref{fig:parabolic/carks}.
	\label{ex:parabolic/subgroups}
\end{example}
\begin{example}[Hyperbolic Subgroups]
	The remaining class of elements in $\psl$ are hyperbolic elements, i.e. those elements which are of absolute trace strictly greater than $2$. Modular graphs corresponding to such elements, named \c{c}ark\index{\c{c}ark} (pronounced chark), again have a unique loop, which is referred to as \emph{spine}. For any hyperbolic element $W$ the Riemann surface $\langle W \rangle \backslash \hh$ is an annulus. There is at least one Farey branch expanding in the direction of the two boundary components, referred to as \emph{inner} and \emph{outer}. The number of Farey branches is determined solely by the element $W$, see Figure~\ref{fig:hyperbolic/carks} for the graph corresponding to the subgroup generated by $W = LSL^{2}S$. 
	\label{ex:hyperbolic/subgroups}
\end{example}
\begin{example}[Dihedral Subgroups]
	A hyperbolic element $W \in \psl$ is called reciprocal if it is conjugate to its own inverse. It turns out that the element conjugating $W$ to $W^{-1}$, denoted $Z_{W}$, is determined uniquely up to multiplication by elements of the group $\langle W \rangle$. $Z_{W}$ is of order $2$. This reflects itself as a special symmetry of the \c{c}ark corresponding to $W$. Namely, if a hyperbolic element $W$ is reciprocal, then the two words, one obtained by tracing the spine of the \c{c}ark clockwise and the other counterclockwise are same, see Figure~\ref{fig:reciprocal} for the graph corresponding to the dihedral subgroup $\langle W,Z_{W}\rangle$. Such subgroups are related to reciprocal geodesics in the modular orbifold which are of some interest from the point of view of analytic number theory, \cite{sarnak/reciprocal/geodesics}. 
\end{example}

\begin{figure}[h!]
	\centering
	\includegraphics[scale=1]{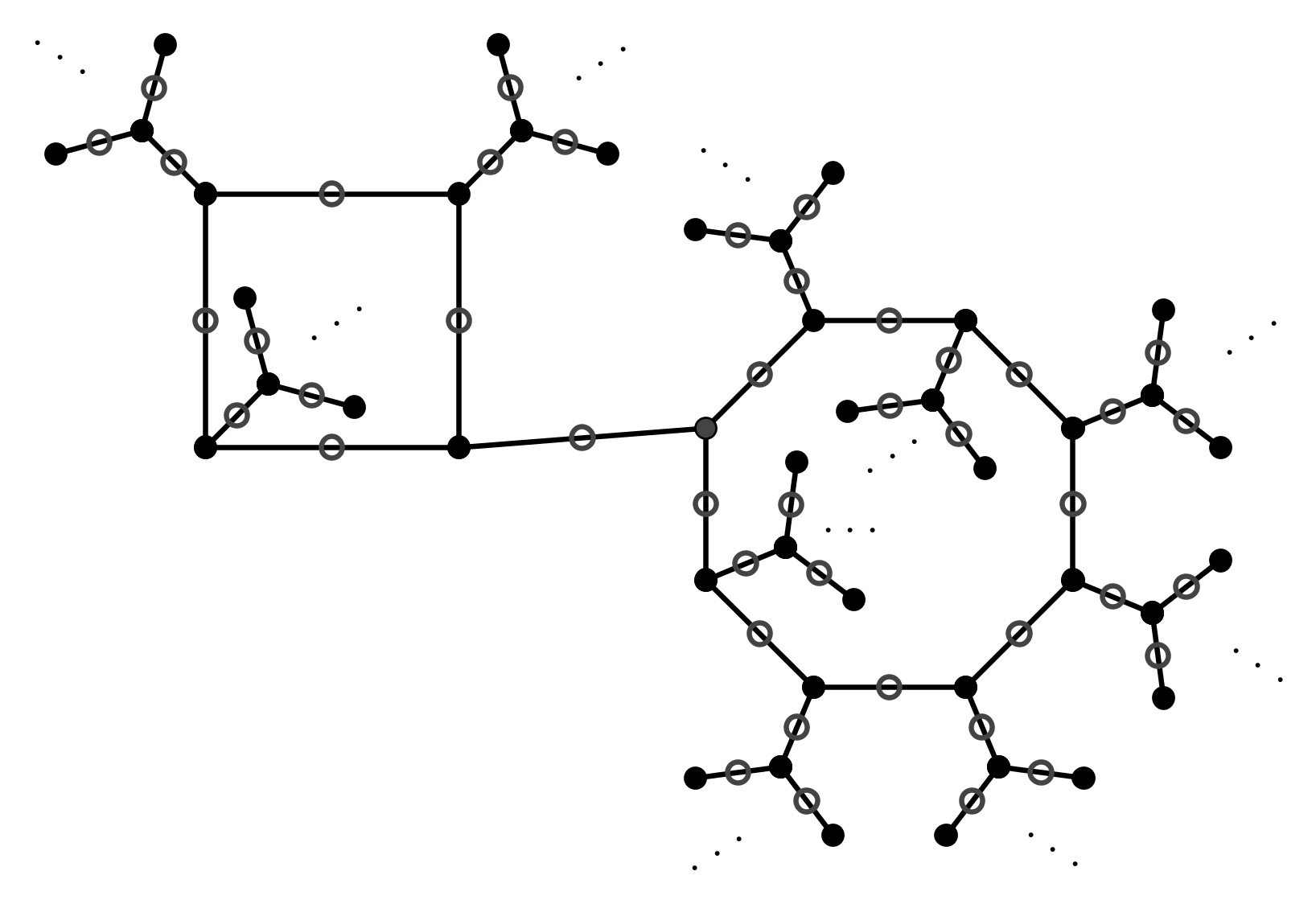}
	\caption{The modular graph corresponding to the pair of pants determined by the subgroup generated by $L^{2}S(LS)^{2}L^{2}SL$ and $ SLSL^{2}SLS(L^{2}S)^{4}(LS)^{2}$.}
	\label{fig:pair/of/pants}
\end{figure}

\begin{example}[Pair of pants]
	Given two hyperbolic elements, say $W_{1}$ and $W_{2}$, which do not commute with each each other, the quotient of the upper half plane by the group $\langle W_{1},W_{2}\rangle$ is a pair of pants which is homeomorphic to a domain in $\CC$ with three boundary components each of which is sa circle. The corresponding modular graph, as in the case of \c{c}arks has Farey tree components expanding in the direction all three boundary components, see Figure~\ref{fig:pair/of/pants} for an example.
\end{example}

\begin{example}[Subgroups isomorphic to the modular group itself]
	Any pair of elements $S'$ and $L'$, of orders $2$ and $3$ respectively, generate a subgroup of $\psl$ which is isomorphic to $\psl$ itself, and this is a proper subgroup whenever the product $S'L'$ is not parabolic. In this case the modular graph 
	$\langle S', L'\rangle\backslash\F$ will be the planar tree with only two terminal vertices, one of type $\circ$ and 
	the other of type $\bullet$.
\end{example}
The above set of examples depict the general picture. That is to say, punctures on the Riemann surface correspond to finite left turn loops on the graph, infinite Farey branches correspond to boundaries of the Riemann surface that are homotopic to the circle, etc.

To summarize the discussion, elements of $\fgicov(\mathcal M)$ are coverings of the modular orbifold by Riemann surfaces, with a finite number of punctures, orbifold points and boundary components. Note that the planar surfaces in $\fgicov(\mathcal M)$ are the Schottky uniformizations of the modular orbifold. Let us finally remark that for any object $\mathcal{X}$ in $\fgicov(\mm)$ the category $\fcov(\mathcal{X})$ consisting of finite covers of $\mathcal{X}$ is in fact a sub-category of $\fgicov(\mm)$. Unlike the ambient category $\fgicov(\mathcal M)$, any two covers in $\fcov(\mathcal X)$ has a common cover in $\fcov(\mathcal X)$.

\section{Arithmetic of modular graphs}
This section is devoted to the study of various arithmetic questions around modular graphs. It turns out that the tame part of the category $\cov(\mm)$ admits an action of $\gal$. The questions here, as addressed in Section~\ref{sec:intro}, are easy to explain. Unexpectedly, there is still very beautiful arithmetic questions around the non-tame part of $\sub(\psl)$ (in fact, concerning the finitely generated infinite index subgroups, ${\rm\bf FGISub}(\psl)$.)

\subsection{Arithmetic of tame subgroups}
The {\it absolute Galois group} $\gal$ of the field of rational numbers is the ultimate product of the classical Galois theory
which is laid before us as a fundamental object of study. This is a non-abelian profinite group,
which is the inverse limit of Galois groups of all finite Galois groups of the field of rationals.
The main open problem of the so-called ``inverse Galois theory'' is to decide whether there exists a
field extension of $\QQ$ with a given finite group as its Galois group. We are far from having
an explicit description of this inverse system of finite groups and we are unable to pin
down an element of the absolute Galois group in an explicit manner, except complex conjugation. Note that being
profinite, the absolute Galois group of $\QQ$ has the structure of a compact, Hausdorff, totally
disconnected topological group. Its size is that of the continuum. 

The way to investigate a complicated group is to construct and study its action on some familiar class of objects.
An action of the group $\gal$ is called a {\it Galois action}\index{Galois action}.
If the action is on a vector space and is continuous, it is usually called a {\it Galois representation}\index{Galois representation}.
A continuous action of a profinite group on a discrete space or a finite-dimensional complex vector space
concerns only a finite quotient of the group and cannot reveal its essential structure if the group is infinite.

In order to study an infinite piece of the absolute Galois group, one idea is to take all
iterates $\{f_i\}_{i=1}^\infty$ of a given polynomial $f$, and study the Galois action on the tree
$T_f$ whose vertices are the points $f_i^{-1}(0)$. In case all the $f_i$s are separable, and $f$ is of degree
$d$, this is a complete rooted $d$-ary tree having a profinite automorphism group ${\rm Aut}(T_f)$,
with a continuous homomorphism $\gal\rightarrow {\rm Aut}(T_f)$. Of course, this representation is
never faithful. A systematic study of this action has started quite recently under the name Arboreal Galois Representations (see \cite{boston/jones/arboreal}). 
We shall not pursue this direction with a dynamical system flavor here, nevertheless trees will enter the scene in a different guise.
Note that the critical points of $f_i$ clutter so the representation we get is in a certain quotient of the fundamental group of the complement of a terrible set.

A technique which proved to be very strong and fruitful to construct and study Galois representations is to use algebraic geometry: $p$-adic \'etale cohomology groups attached to an algebraic variety are finite-dimensional $\QQ_p$-vector spaces with a continuous action of $\gal$, and the image of such representations may be infinite, although ``small" in the sense that
such representations always have a ``big" kernel. The main open problem of this field is the characterization of representations of geometric origin among all linear representations (Fontaine-Mazur conjecture
\cite{kisin/galois/representation}).

A natural next step is to study the Galois actions on (quotients of) algebraic fundamental groups. These are called ``large" representations because they may be faithful. \'Etale fundamental groups were introduced and studied by Grothendieck in the late sixties \cite{rev-etale}. In this setup, if $k$ is a field, then the (\'{e}tale) fundamental group of ${\rm Spec}(k)$ is $\mathrm{Gal}(k)$. In the arithmetic case one cannot employ the usual topological tools to define the fundamental group. However, there is a good notion of covering, i.e. \'etale covering, and one uses the category of all \'etale coverings  of the space in question, in order to define its fundamental group.

\subsubsection{Galois action on algebraic fundamental groups.} \label{sec:Galoisactionsonpi1} Let
$\overline{\QQ}\subset \CC$ be the algebraic closure of $\QQ$ in $\CC$ and let $X_\QQ$ be a smooth
geometrically irreducible variety over $\QQ$ (it is more natural to consider $X$ to be a stack). 
Let $p\in X$ be a rational point (in the scheme-theoretic sense) and let
$\overline{p}$ be a geometric point of $X_{\overline{\QQ}}:=X\!\otimes_{_\QQ}\!\! \overline{\QQ}$ above
$p$. One associates \'etale fundamental groups to $(X,p)$ and to $(X_{\overline{\QQ}}, \overline{p})$ together
with an exact sequence
\begin{equation}\label{exact/1}
1\rightarrow \pi_1(X_{\overline{\QQ}},\overline{p}) \rightarrow \pi_1(X_\QQ,p)
\rightarrow \gal \rightarrow 1
\end{equation}
which defines an {\it outer action} $\rho_{X,p}: \gal \rightarrow
\mathrm{Out}(\pi_1(X_{\overline{\QQ}},\overline{p}))$. The group $\pi_1(X_{\overline{\QQ}},\overline{p})$ is canonically
isomorphic to the profinite completion of the topological fundamental group of $X_\CC$, that is,
$\pi_1(X_{\overline{\QQ}},\overline{p})=\hat{\pi}_{1, top}(X_\CC,p)$.
Here $X_\CC$ denotes the complexification of $X$. 

The simplest interesting situation is $X=\PP^1_\QQ \backslash\{0,\infty\}$ 
($=:\mbox{Spec}\, \QQ[x,1/x]$). In this case the exact sequence is  
\begin{equation}
1\rightarrow \pi_1(\PP^1_{\overline{\QQ}}\backslash\{0,\infty\},\overline{p}) \rightarrow \pi_1(\PP^1_\QQ\backslash\{0,\infty\},p).
\rightarrow \gal \rightarrow 1
\end{equation} 
Since the topological fundamental group of $\PP^1_\CC \backslash\{0,\infty\}$ is $\ZZ$, 
the algebraic fundamental group $\pi_1(\PP^1_{\overline{\QQ}},\overline{p})$ 
is isomorphic to the profinite completion of integers $\widehat{\ZZ}$. 
In fact, we can describe the system of all coverings of $X$ explicitly, by the power maps 
$\varphi_n:X\stackrel{z^n}{\rightarrow}X$.
The Galois group acts on the monodromy group $\ZZ/n\ZZ$ of the covering $\varphi_n$ through the cyclotomic character $\sigma\in \gal\rightsquigarrow \chi(\sigma)\in \widehat{\ZZ}$, such that 
$\chi(\sigma):[i]\in \ZZ/n\ZZ \rightarrow [i\chi(\sigma)]\in \ZZ/n\ZZ$. (If we represent $\ZZ/n\ZZ$ as the group of $n$th roots of unity in $\CC$ generated by a primitive $n$th root $\xi$, 
the action is described by $\xi^i\rightarrow \xi^{i\chi(\sigma)}$). For a given element $\sigma$, 
this defines an automorphism of each monodromy group, in a way compatible  with the system of coverings $\{\phi_n\}$. 
To sum up, the above exact sequence becomes 
\begin{equation}
1\rightarrow \widehat{\ZZ} \rightarrow \pi_1(\PP^1_\QQ \backslash\{0,\infty\})
\rightarrow \gal \rightarrow 1
\end{equation}
inducing a surjective homomorphism $\gal \rightarrow \mathrm{Aut}(\widehat{\ZZ})\simeq \widehat{\ZZ}^\times$.
In fact, this is the abelianization map of $\gal$. 

The simplest non-abelian example of a Galois action on a fundamental group
is the case of the dihedral triangle orbifold with signatures $(2,2,\infty)$. 
However, this requires some preparation to orbifolds (or stacks), so let's make the situation non-abelian by removing one more point from $\PP^1_\QQ \backslash\{0,\infty\}$.
So let $Y=\PP^1_\QQ \backslash\{0,1,\infty\}$ ($=:\mbox{Spec}\, \QQ[x,,1-x, 1/x)]$).
Since the topological fundamental group of $\PP^1_\CC \backslash\{0,1,\infty\}$ is free of rank 2, the group $\pi_1(\PP^1_{\overline{\QQ}}\backslash\{0,1,\infty\},\overline{p})$ is the free profinite group $\widehat{\F}_2$ of rank 2. 
To sum up, the above exact sequence becomes 
\begin{equation}\label{exact/2}
1\rightarrow \widehat{\F}_2 \rightarrow \pi_1(\PP^1_\QQ \backslash\{0,1,\infty\})
\rightarrow \gal \rightarrow 1.
\end{equation}
So any element of $\gal$ determines an automorphism of $\widehat{\F}_2$, up to an inner automorphism.
There is a standard way to lift the {\it outer} Galois action to an {\it outer} action $\mathrm{Aut}(\pi_1(\overline{X},\overline{p}))\simeq \mathrm{Aut}(\widehat{\F}_2)$.
In other words, there is a homomorphism $\gal \rightarrow \mathrm{Out}(\widehat{\F}_2)$.

Belyi's theorem implies that this representation is faithful, and the characterization of its image is the main question of Grothendieck-Teichm\"uller
theory (GT) initiated by Belyi, Drinfeld and Ihara (see \cite{schnepsgt} and \cite{lochak/fragments}). As a result of decades of intense research by prominent mathematicians, several sets of equations that must be respected by this image have been found. At the current level of research, it seems an intractable question to see whether these equations are independent or whether they are sufficient to characterize the image. 

We shall probably never have a satisfactory understanding of the group $\mathrm{Aut}(\widehat{\F}_2$) and the associated Galois action.
One may say that the Galois representations on fundamental groups
are ``too large". A remedy to this problem is to consider representations on smaller quotients of
fundamental groups, such as the pro-nilpotent or pro-l fundamental
groups. At this point starts the motivic side of the theory, initiated by Wojtkowiak and Deligne
(see \cite{deligne/goncharov/motives} and \cite{hain/weighted/completion/of/galois}). Although this is a very rich territory, with connections to
multiple zeta values, polylogarithms, etc. it is still ``linear" and ``far from the anabelian dream
imagined by Grothendieck." For a discussion of the ``non-linear" theory, see \cite{lochak/fragments}. 
The hypergeometric Galois actions that we proposed in the previous chapter might provide a
possibility to go beyond the linear theory.

\subsubsection{Galois action on modular graphs.}\label{sec:galois/action/on/graphs}

If $G$ is a finite index subgroup of $\psl $, then the projection $\pi_{G}$ is a Belyi map, i.e. it is a covering map from the Riemann surface $G \backslash \hh$ to the modular orbifold $\mathcal M$. The latter is conformally equivalent to $\mathbb C$. 
Compactifying by adding a point at infinity, we get the Riemann sphere and  the orbifold covering
$$
\pi_{G} \colon G \backslash \hh \to \mathcal M
$$ 
can be viewed as a branched covering 
$$
\overline{\pi_{G}}: \overline{G \backslash \hh} \to \mathbb P^1(\CC)
$$ 
branched at elliptic points and at infinity. Branching at elliptic points is restricted whereas there is no restriction on the branching above the point at infinity. This is no loss of generality as there is a degree $6$ covering from $\PP^{1} \sm \{0,1,\infty\}$ to the modular orbifold. By Belyi's theorem this branched covering admits a model whose defining polynomials have coefficients from a number field (i.e. a finite extension of $\QQ$)\footnote{A converse to this claim was first stated by Belyi whose \emph{proof} relied entirely on the article by Weil, \cite{weil/field/of/definition}. Experts in the area believed that some explanation was necessary, see \cite{juergen/obvious/part/of/belyi}}. 
In particular, $G \backslash \hh$ is defined over a number field. Such curves will be called \emph{arithmetic curves}. This allows us to define an action of the absolute Galois group on $\mathbf{FCov}(\mathcal M)$, or equivalently on the category of modular graphs.

By Belyi's theorem every arithmetic curve arises as a covering of \linebreak $\PP^{1} \sm \{0,1,\infty\}$ and hence the modular orbifold: arithmetic curves are precisely the compactifications of the covering curves inside $\mathbf{FCov}(\mathcal M)$. As we already noted above, every arithmetic curve will appear in $\mathbf{FCov}(\mathcal M)$ in infinitely many ways, i.e. $\mathbb P^1(\CC)$ will appear as many times as the number of trivalent graphs drawn on the sphere. The last example illustrates why we consider the action of the Galois group on the coverings (curve plus the covering), and not just on the curves. The arithmetic curve here is just $\mathbb P^1(\CC)$ which is defined over $\QQ$ and so the action is trivial, whereas the action becomes faithful when we consider the coverings as well.

The action of $\gal$ preserves (see \cite{jones/streit/galois/groups}):

\begin{enumerate}
	\item the number of edges in the modular graph, equivalently, the degree of the branched covering;
	\item the genus of the arithmetic curve;
	\item the degree distribution of vertices of type $\circ$, equivalently, the branching behavior above $\circ$ (recall that this degree can be only 1 or 2);
	\item the degree distribution of vertices of type $\bullet$, equivalently, the branching behavior above $\bullet$ (recall that this degree can be only 1 or 3);
	\item the degree distribution of right-turn cycles, equivalently, the branching behavior above the point at infinity (there is no restriction on this);
	\item monodromy group.
\end{enumerate}

As the aim is to understand $\gal$, one looks for families of arithmetic curves supporting a faithful action of $\gal$. It appears that there are many candidates for this. Motivated by the above list one may first try fixing the genus, $g$. In fact, for any non-trivial element $\sigma$ of $\gal$ one can easily find an algebraic number, say $\alpha$, on which $\sigma$ acts non-trivially. Then, if $E$ is an elliptic curve with $j$ invariant $\alpha$, then $E$ and $E^{\sigma}$ are distinct elliptic curves. Hence genus $1$ arithmetic curves admit a faithful $\gal$ action. By Belyi's theorem every such curve corresponds to a modular graph of genus $1$. The conclusion is that the Galois action on the set of modular graphs of genus 1 is faithful. One can also show that the same is true for genus $0$ curves. More generally, we have:
\begin{theorem}[{\cite{ernest/gabino/genus/g/Galois/action}}]
	Fix a non-negative integer $g$. The action of the absolute Galois group on the set of genus $g$ modular graphs is faithful. 
\end{theorem}

In a different vein, the Galois action sends planar trees to planar trees by the above conservation properties, and we have the result:
\begin{theorem}[\cite{sch-dd-on-rsphere}]
	The action of the absolute Galois group on the set of modular graphs which are planar trees is faithful.
\end{theorem}

There are a myriad of  questions  related to the arithmetic surrounding the theory built so far,  and there is a vast literature about these questions. We refer the interested reader to \cite{ernesto/gabino/book} for further reading. We do not pursue here this ``{finite}" side of the theory further, which is devoted to the study of $\fcov(\mathcal M)$ of {\it finite} coverings of the modular orbifold, and we turn our attention to the study of some special infinite covers.

\subsection{Arithmetic of non-tame subgroups}

At this stage, the reader will probably expect us to define a Galois action on $\mathbf{FCov}(\mathcal X)$, where $\mathcal X$ is an object in $\cov(\mm)$ exactly as in the finite case. Likewise for the ambient category $\mathbf{FGICov}(\mathcal M)$. This requires an interpretation of the elements and morphisms of this category as arithmetic objects over some field, which is possibly some transcendental extension of $\QQ$ depending on $\mathcal X$.  We firmly believe that such an interpretation does exist and is fruitful, however, we don't know how to make this. Instead we will give an alternative arithmetic interpretation of hyperbolic $\ZZ$-covers from \cite{UZD}.

\subsubsection{Hyperbolic subgroups and binary quadratic forms.} 

Let $W = \mat{p&q}{r&s}$ be a hyperbolic element in $\psl$, i.e. $|p+s|>2$. The action of $W$ on $\HH$ has two fixed points, determined by roots of the following integral equation:
$$rx^{2} + (s-p)x -q = 0.$$
Since solutions of the above equation remains fixed when one takes non-zero multiples of the coefficients, the above equation can be reduced to:
\begin{equation}
	\frac{r}{\delta}x^{2} + \frac{s-p}{\delta}x - \frac{q}{\delta} = 0,
	\label{eq:roots}
\end{equation}
where $\delta$ is the greatest common divisor of $r$, $s-p$ and $q$. Note that hyperbolicity of the matrix implies that the equation has two real roots. There is a unique geodesic in $\HH$ joining these two fixed points, called the \emph{geodesic of} $W$. Indeed, the action of $W$ on $\HH$ is nothing but a \emph{translation} along this geodesic. 

The projectivization of the left hand side of Equation~\ref{eq:roots}, namely
$$\frac{r}{\delta}x^{2} + \frac{s-p}{\delta}xy - \frac{q}{\delta}y^{2}$$
is then an indefinite binary quadratic form\index{indefinite binary quadratic form}\index{binary quadratic form!indefinite}, i.e. a homogeneous of degree two element of $\ZZ[x,y]$ of positive discriminant. This form will be denoted by $f_{W}$. Conversely, given any triple $(a,b,c)$ with $b^{2} - 4ac$ being positive and the greatest common divisor of $a$, $b$ and $c$ being $1$, by solving a suitable Pell equation, one can find a hyperbolic element of $\psl$, $W$, so that 
$$f_{W} = ax^{2} + bxy + cy^{2}.$$

On the other hand, $\psl$ acts on the set of indefinite binary quadratic forms by change of variable. This action transforms into the conjugation action in the language of subgroups. We have:

\begin{theorem}[{\cite[Theorem 3.1]{UZD}}]
	There is a one to one correspondence between $\psl$-classes of binary quadratic forms and conjugacy classes of subgroups of $\psl$ generated by one hyperbolic element.
	\label{thm:correspondence/between/forms/and/subgroups}
\end{theorem}

As described in Example~\ref{ex:hyperbolic/subgroups} such modular graphs have a unique loop (called the spine), and finitely many Farey branches expanding in the direction of both connected components of the boundary of the annulus, $\langle W \rangle \backslash \HH$. 

Given two forms $f_{1} = a_{1}x^{2} + b_{1}xy + c_{1}y^{2}$  and $f_{2} = a_{2}x^{2} + b_{2}xy + c_{2}y^{2}$ there is a product defined by Gau{\ss}, \cite{disquisitiones}, if and only if either $\Delta_{1}/\Delta_{2} = (b_{1}^{2} - 4a_{1}c_{1})/(b_{2}^{2} - 4a_{2}c_{2})$ or $\Delta_{2}/\Delta_{1}$  is a perfect square. In particular the set of all forms of discriminant $\Delta$ is a group under this operation. The function mapping a binary quadratic form $f = ax^{2} + bxy + cy^{2}$ of discriminant $\Delta = b^{2} - 4ac$ to the narrow ideal class generated by $(1,\omega)$, where $\omega = \frac{b + \sqrt{\Delta}}{2a}$ in the quadratic number field $\QQ(\sqrt{\Delta})$ is a group isomorphism. Therefore, given a square-free positive integer $\Delta$, finding the class number of the number field $\QQ(\sqrt{\Delta})$, i.e. the size of the class group of $\QQ(\sqrt{\Delta})$, is equivalent to finding the number of inequivalent $\psl$ classes of indefinite binary quadratic forms of same discriminant. The latter set of problems, which are now older than 200 years, the class number problems of Gau{\ss}, and are all stated in Disquisitiones Arithmeticae, \cite{disquisitiones}. We must mention that Gau{\ss} was not able to solve the analogous problems for negative discriminant but guessed almost all of class number one. 

\begin{figure}[h!]
	\centering
	\includegraphics[width=12cm]{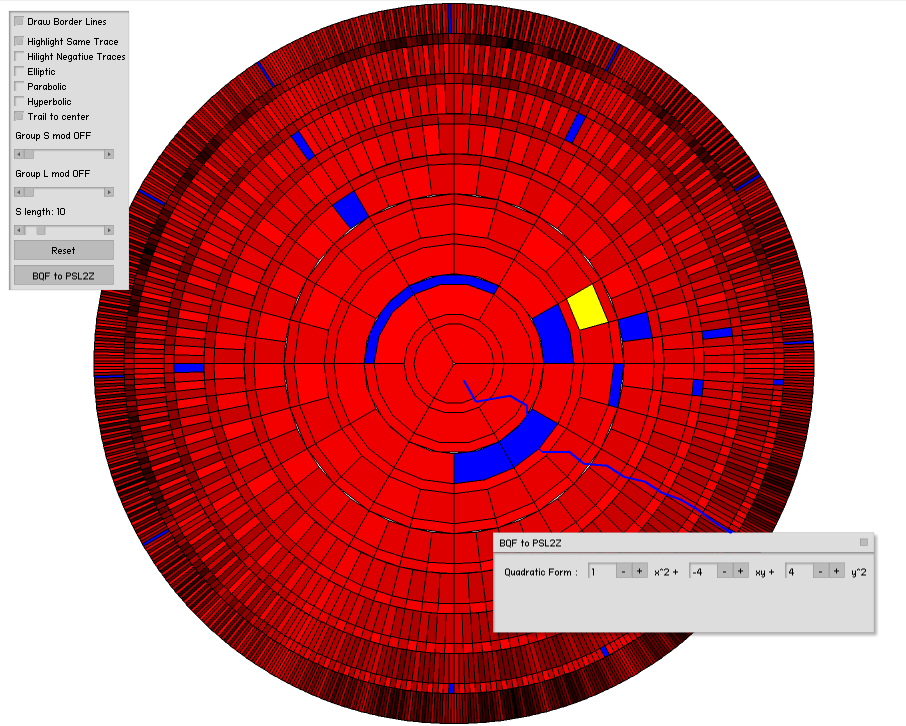}
	\caption{The application InfoMod}
	\label{fig:infomod}
\end{figure}

There are further arithmetic problems concerning indefinite binary quadratic forms:

\begin{enumerate}
	\item minimum problem: find the smallest positive integer $N$ which is attained, and
	\item representation problem: given an integer $N$ determine whether the equation $f = ax^{2} + bxy + cy^{2} = N$ has ingtegral solutions, and if exists how many.
\end{enumerate}

\noindent Both of the solutions have beautiful interpretations and algorithmic solutions in terms of \c{c}arks, see \cite{reduction}. Moreover the corresponding geodesics and above arithmetic problems are visualized in the application InfoMod, see Figure~\ref{fig:infomod}, by the authors and H. Ayral. 


\medskip
\noindent {\bf Acknowledgements.} 
We are thankful to Athanase Papadopoulos for inviting us to publish in this volume and for his comments on the previous versions of this manuscript.

\bibliographystyle{abbrv}

\end{document}